\documentclass[letterpaper, 10 pt, conference]{ieeeconf}  

\IEEEoverridecommandlockouts                              

\usepackage{psfrag}
\usepackage{graphicx}
\usepackage{cite}
\usepackage{amsmath} 
\usepackage{amssymb}  
\usepackage{amsfonts}
\usepackage{mathrsfs}
\usepackage{mathtools}
\usepackage[dvipsnames]{xcolor}
\usepackage[thmmarks, amsmath]{ntheorem}
\usepackage{circuitikz}
\usepackage{booktabs}
\usepackage{colortbl}
\usepackage{pgfplots}
\usetikzlibrary{arrows}
\usepackage{marvosym}
\usepackage[ruled,vlined]{algorithm2e}
\usepackage{enumerate}

\renewcommand{\baselinestretch}{0.99}

\usepackage{soul}

\newcommand{\setdef}[2]{\{#1 \;|\; #2\}}

\usepackage{tikz}
\usepackage{pgfplots}
\usetikzlibrary{arrows}
\usetikzlibrary{decorations.pathmorphing}
\usetikzlibrary{decorations.markings}
\usetikzlibrary{snakes}
\usetikzlibrary{matrix}
\usetikzlibrary{calc}
\usetikzlibrary{patterns}
\usetikzlibrary{positioning}
\usetikzlibrary{shapes}
\usetikzlibrary{spy}
\usetikzlibrary{shapes.symbols}
\usetikzlibrary{shapes.geometric}
\usetikzlibrary{shadows.blur}
\tikzstyle{every node}=[font=\small]
\tikzstyle{every path}=[line width=0.8pt,line cap=round,line join=round]

\usepackage{ifthen,calc}
\usepackage{xstring}

\usetikzlibrary {positioning}
\usepackage{empheq}
\usepackage[most]{tcolorbox}

\newtcbox{\mybox}[1][]{%
    nobeforeafter, math upper, tcbox raise base,
    enhanced, colframe=blue!30!black,
    colback=blue!30, boxrule=1pt,
    #1}

\pgfplotsset{compat=1.14}

\newcommand{\vect}[1]{\mathbbold{#1}}

\newcommand{\vzeros}[1][]{\vect{0}_{#1}}
\DeclareSymbolFont{bbold}{U}{bbold}{m}{n}
\DeclareSymbolFontAlphabet{\mathbbold}{bbold}

\DeclareMathOperator*{\argmin}{argmin}

\providecommand{\norm}[1]{\lVert#1\rVert}
\newcommand{\R}{{\mathbb R}}
\newcommand{\mc}{\mathcal}

\newtheorem{proposition}{Proposition}
\newtheorem{assumption}{Assumption}
\newtheorem{remark}{Remark}
\newtheorem{definition}{Definition}

\graphicspath{{Figures/} }

\newcommand{\tr}{^{\sf T}}

\newcommand{\vi}[2]{\mathrm{VI}(#1,#2)}

\usepackage{scalerel,stackengine}
\stackMath
\newcommand\reallywidehat[1]{%
\savestack{\tmpbox}{\stretchto{%
  \scaleto{%
    \scalerel*[\widthof{\ensuremath{#1}}]{\kern.1pt\mathchar"0362\kern.1pt}%
    {\rule{0ex}{\textheight}}
  }{\textheight}%
}{2.4ex}}%
\stackon[-6.9pt]{#1}{\tmpbox}%
}

\SetKwBlock{Iterate}{Iterate}{End}

\definecolor{pinegreen}{cmyk}{0.92,0,0.59,0.25}
\definecolor{royalblue}{cmyk}{1,0.50,0,0}
\definecolor{lavander}{cmyk}{0,0.48,0,0}
\definecolor{violet}{cmyk}{0.79,0.88,0,0}
\tikzstyle{cgray}=[circle, draw, thin,fill=black!20, scale=0.8]
\tikzstyle{qblack}=[rectangle, draw, thin,fill=black!90, scale=0.8]
\tikzstyle{gpath}=[ultra thick, draw, black!70, opacity=1, ->]                       
\tikzstyle{legend_overlay}=[rectangle, rounded corners, thin,
                           top color= white,bottom color=green!25,
                           minimum width=2.5cm, minimum height=0.8cm,
                           pinegreen]
\tikzstyle{legend_phytop}=[rectangle, rounded corners, thin,
                          top color= white,bottom color=cyan!25,
                          minimum width=2.5cm, minimum height=0.8cm,
                          royalblue]
\tikzstyle{legend_general}=[rectangle, rounded corners, thin,
                          top color= white,bottom color=lavander!25,
                          minimum width=2.5cm, minimum height=0.8cm,
                          violet]
                          
\tikzstyle{block_m}=[draw, fill=black!20, text width=8em, 
    text centered, minimum height=3em,drop shadow]

\tikzstyle{block_s}=[draw, fill=black!20, text width=10em, 
    text centered, minimum height=3em,drop shadow]    

\begin{document}

\title{\LARGE \bf Towards robustness guarantees for feedback-based optimization
\thanks{
M. Colombino and A. Bernstein are with the National Renewable Energy Laboratory, Golden, Colorado, \{{name.lastname}\}@nrel.gov; J. W. Simpson-Porco is with the Department of Electrical and Computer Engineering at the University of Waterloo, jwsimpson@uwaterloo.ca; This work was authored in part by the National Renewable Energy Laboratory, operated by Alliance for Sustainable Energy, LLC, for the U.S. Department of Energy (DOE) under Contract No. DE-AC36-08GO28308. Funding provided by DOE Office of Electricity, Advanced Grid Modeling Program, through agreement NO. 33652. The views expressed in the article do not necessarily represent the views of the DOE or the U.S. Government. The publisher, by accepting the article for publication, acknowledges that the U.S. Government retains a nonexclusive, paid-up, irrevocable, worldwide license to publish or reproduce the published form of this work, or allow others to do so, for U.S. Government purposes.
} }
\author{Marcello Colombino, John W. Simpson-Porco, Andrey Bernstein}

\thispagestyle{empty}
\pagestyle{empty}

\maketitle

\begin{abstract}
Feedback-based online optimization algorithms have gained traction in recent years because of their simple implementation, their ability to reject disturbances in real time, and their increased robustness to model mismatch. While the robustness properties have been observed both in simulation and experimental results, the theoretical analysis in the literature is mostly limited to nominal conditions. In this work, we propose a framework to systematically assess the robust stability of feedback-based online optimization algorithms. We leverage tools from monotone operator theory, variational inequalities and classical robust control to obtain tractable numerical tests that guarantee robust convergence properties of online algorithms in feedback with a physical system, even in the presence of disturbances and model uncertainty. {The results are illustrated via an academic example and a case study of a power distribution system.}
\end{abstract}

\section{Introduction}

%
Online optimization methods are traditionally well suited for classical computer science tasks (recommendation engines, classifications. etc.) but are not usually designed to deal with the complex constraints and safety requirements of physical systems affected by unknown disturbances. A prototypical example of such system is the electric power grid, where the system operator must optimally schedule power generation while taking into account capacity constraints, voltage/current safety constraints, and unknown disturbances represented by uncontrollable loads and variable generation. Traditionally, the optimal operation of such large scale engineering systems is done via frequent re-optimization based on complex models and disturbance forecasts. Recently, however, much simpler online (or feedback-based) optimization methods have been proposed for constrained engineering systems with tremendous success in applications ranging from communication networks \cite{low1999optimization} to power systems~\cite{li2016connecting,Bolognani_feedback_15,Dallanese2016optimal,bernstein2019online,tang2017real,hauswirth2016projected} to transportation~\cite{MV-JC:18-necsys}. 

The appeal of using feedback-based over off-line optimization approaches is the same appeal of using feedback over feedforward control: feedback optimization methods show superior robustness to model uncertainty and are able to attenuate or reject unmeasured disturbances. While these properties have been extensively observed in the literature, the analysis of these algorithms is mostly performed under nominal conditions. Some exceptions are~\cite{colombino2018online,lawrence2018optimal}, where robustness to linear unmodeled dynamics is considered.

In this work, we focus on a simple first-order online approximate gradient descent, similar to those proposed in~\cite{hauswirth2016projected,Dallanese2016optimal}. The novel contributions of the work are that we
\begin{itemize}
\item characterize the equilibria of the feedback interconnection of the physical system and the online optimization scheme using methods from the literature on variational inequalities~\cite{facchinei2007finite} and monotone operators~\cite{bauschke2011convex};
\item propose a framework, based on classical robust control theory~\cite{GED-FP:00,AM-AR:97,CS-SW:15,JV-CWS-HK:16} that allows us to systematically test the robustness properties of the online algorithm by guaranteeing robust stability with respect to a large class of uncertain physical systems;
\item validate our results both on an academic example and a case study of a power distribution system, for which the proposed methods allows us to verify robust stability for a wide range of realistic operating conditions.
\end{itemize}
The robustness analysis in this work is partly inspired by~\cite{lessard2016analysis}, even though we consider uncertainty in the model of a physical system and not merely as a tool to characterize the properties of certain nonlinear operators. The analysis in this paper differs from~\cite{Dallanese2016optimal} as it is not based on a bounded error between the output of the real nonlinear system and the one of an approximate linear system. This leads to sharper convergence guarantees (to a point instead of a set). Similarly to~\cite{Dallanese2016optimal,bernstein2019online} all results in this paper carry over to the time-varying setting with minimal modifications. 

\emph{Notation: }
For a symmetric positive definite $P \succ 0$,  $\langle \cdot, \cdot \rangle_{P}: \R^n \times \R^n \rightarrow \R$ denotes the inner product $\langle x, y \rangle_{P} := x\tr P y$, and the corresponding induced norm $\left\| \cdot \right\|_P : \R^n \rightarrow \R_{\geq 0}$ defined as $\left\| x \right\|_P := \sqrt{ x\tr P x }$. 
A map $f: \R^n \rightarrow \R^n$ is $L$-Lipschitz w.r.t the inner product $\langle \cdot,\cdot\rangle_{P}$, if $\exists\,L>0$ such that $\left\| f(x) - f(y)\right\|_P \leq L \left\| x-y\right\|_P$ for all $x,y \in \R^n$. If $f$ is differentiable, then $\partial f(x)$ denotes its Jacobian matrix at $x$. A function $f:\R^n\to \R$ is $L$-strongly smooth w.r.t $\langle \cdot,\cdot\rangle_{P}$, if it is differentiable and its gradient is $L$-Lipschitz w.r.t $\langle \cdot,\cdot\rangle_{P}$. Every mentioned set $\mathcal{S} \subseteq \R^n$ is nonempty. For a closed convex set $\mathcal{C} \subseteq \R^n$, the projection operator, $\text{Proj}_{\mathcal{C}}^{P} : \R^n \rightarrow \mathcal{C} \subseteq \R^n $, is defined as $\text{Proj}_{\mathcal{C}}^{P}(x) := \arg \min_{y \in \mathcal{C}} \left\|x-y \right\|_P$.  If $P$ is omitted it is assumed that $P=I$.

\section{Preliminaries}

This section introduces some preliminary results from the theories of monotone operators and variational inequalities (VIs), which are central in many fields of applied mathematics, engineering and economics, We refer the reader to~\cite[Chapter 1]{gentile2018equilibria} for a gentle introduction and~\cite{bauschke2011convex,facchinei2007finite} for a comprehensive review.
%
%
We begin by defining a monotone operator.
\begin{definition}[{Monotone operator}]\label{def:MON}
Given $P\succ 0 $ and $\rho > 0$, a map $F: \mc S\subset\R^n \rightarrow \R^n$ is $\rho$-strongly monotone on $\mathcal{S}$ w.r.t $\langle \cdot,\cdot\rangle_{P}$ if, $\forall\,x, y \in \mc S,$ 
$
\langle  x-y,F(x) - F(y)\rangle_{P} \geq \rho\norm{x-y}_P^2,
$
We say $F$ is monotone if the inequality holds for $\rho=0$.
\end{definition}
%
Next we define the $\emph{Clarke generalized Jacobian}$ for non-smooth locally Lipschitz maps.
\begin{definition}[{Clarke generalized Jacobian}~\cite{clarke1990optimization}] \label{def:CJ}
Let ${\mathcal S \subseteq \R^n}$ be a closed-convex set and  $F:\mathcal S\to \R^n$ be a locally Lipschitz map. The \emph{Clarke generalized Jacobian} of $F$ at $x \in \mathcal{S}$ is defined as the set 
\[
\partial_C F(x) =  \mathrm{co}\Bigg\{ J \in \R^{n \times n} :  J = \lim_{
\substack{x_i\to x \\F(x_i) \textnormal{ differentiable}}} \partial F(x_i)\Bigg\}
\]
where $\mathrm{co}$ is the convex hull. { For a continuously differentiable function $g:\mathcal{S}\to \R$, its Clarke generalized Hessian at $x \in \mathcal{S}$ is defined as the set $\partial^2_C g(x) = \partial_C\nabla g(x).$ }
\end{definition}
For non-differentiable and locally Lipschitz maps, the following proposition characterizes monotonicity in terms of the Clarke generalized Jacobian.

\begin{proposition}[Characterizing monotonicity]\label{propo.jen.jac}
Let $\mathcal S \subseteq \R^n$ be a closed-convex set and  $F:\mathcal S\to \R^n$ be a locally Lipschitz map. Then $F$ is $\rho$-strongly monotone {on $\mathcal{S}$} w.r.t $\langle \cdot,\cdot\rangle_{P}$ if and only if
\begin{equation} \label{eq:SMON-inequality.propo.jen}
\frac{1}{2}\left[J\tr P + P\,J \right]\succcurlyeq \rho P
\end{equation}
for all $J\in\partial_CF(x)$ and all $x \in \mc S$, and is monotone iff~\eqref{eq:SMON-inequality.propo.jen} holds for $\rho=0$.
\end{proposition}

Proposition~\ref{propo.jen.jac} is a slight generalization to~\cite[Proposition 2.1]{schaible1996generalized} for monotonicity w.r.t $\langle \cdot,\cdot\rangle_{P}$.
Next, we recall the idea of a VI, which is an important tool for the rest of the paper
\begin{definition}[VI Solution]
Consider a set $\mathcal X \subseteq \R^n$ and a map $F:\mathcal X\to \R^n$. A point $\bar x\in\R^n$ is a solution of the variational inequality $\vi{\mc X}{F}$ w.r.t $\langle \cdot,\cdot\rangle_{P}$ if 
\[
\langle F(\bar x),x-\bar x\rangle_{P} \ge 0,\quad \forall\,x\in\mc X\,.
\]
\end{definition}
The following algorithm can be used to find the solution of VIs involving a strongly monotone operator $F$.
\begin{proposition}[Convergence of the projection algorithm]{\cite[Theorem 12.1.2]{facchinei2007finite}}\label{propo.conv.grad} Let $\mc X \subseteq \R^n$ be closed and convex, and let $F:\mc X\to \R^n$ be $\rho$-strongly monotone and $L$-Lipschitz on $\mathcal{X}$ w.r.t $\langle \cdot,\cdot\rangle_{P}$. Then the variational inequality $\vi{\mc X}{F}$ admits a unique solution $x^\star$, and, for any { $x_1 \in \mathcal{X}$}, the sequence $\{x_k\}_{k=1}^\infty$ generated by the projection algorithm
\[
x_{k+1} =  \mathrm{Proj}^P_{\mc X}\left(x_k -\tau F(x_k) \right)\,, \quad x_1 \in \mathcal{X}\,.
\]
 with $\tau < \frac{2\rho}{L^2}$ converges { geometrically} to $x^\star$, Moreover, the best geometric convergence rate of $1-(\rho/L)^2$ is achieved with step size $\tau = \rho/L^2$.
\end{proposition}
\section{Feedforward vs feedback optimization}\label{sec.setup}
In this work, we consider a physical system that maps an input $u\in\mc U \subseteq \R^n$ and an unknown disturbance $w\in\mc W\subseteq \R^p$ to the output $y\in\R^m$ according to the map
\begin{align}\label{eq.physical.system}
y =  \pi(u,w),
\end{align}
where $\pi:\mc U\times \mc W\to\R^m$ is continuously differentiable and {locally} Lipschitz continuous in $u$. A system operator is responsible for managing the physical system~\eqref{eq.physical.system}, and is interested in solving a optimization problem of the form
\begin{align}\label{eq.opt.orig}
\begin{split}
    \min_{u\in\mc U} ~&~f(u)\\
    \mathrm{subject~to}&~y = \pi(u,w),~ y\in\mc Y
    \end{split}
\end{align}
to achieve optimal operation of~\eqref{eq.physical.system}. In~\eqref{eq.opt.orig} $u\in\R^n$ is the decision variable that we apply to the system, $w\in\mc W$ is an unknown disturbance, $\mc U\subset \R^n$ is a {closed} convex set that represents hard physical limits on the input $u$, $f:\mc U\to \R$ is a strongly smooth convex function, and $\mc Y$ represents (for example) safety constraints on the system's output $y$. 
\subsection{The standard approach: feedforward optimization} \label{sec.feedforward}
Optimization problems of the form~\eqref{eq.opt.orig} are ubiquitous in many engineering disciplines (e.g., optimal power or gas flow, optimal traffic control, etc.) and are in general hard to solve. This is because they are often of very large size, non-convex ($\pi$ is generally nonlinear) and they require a precise knowledge of $\pi$ and of the disturbance $w$, which is often unavailable. In many real-world applications, the system operator will have access to a linearized model of the system \eqref{eq.physical.system} of the form 
\[
y \approx { \Pi} u + \Pi_w w
\]
and to a forecast or guess $\hat w$ for the disturbance $w$. The operator can then periodically solve the \emph{feedforward} optimization problem\footnote{Or perhaps a robust or chance-constrained version to account for model mismatch and forecast error.}
\begin{align}\label{eq.opt}
\begin{split}
 \hat u =   \argmin_{u\in\mc U} ~&~f(u)\\
    \mathrm{subject~to}&~y =  { \Pi} u + \Pi_w \hat w,~ y\in\mc Y
    \end{split}
\end{align}
in order to protect themselves from model uncertainty and forecast error. Solving~\eqref{eq.opt} and applying the solution $\hat u$ to the system can be seen as analogous to applying feedforward control to a dynamical system. The solution $\hat u$ is the best solution we obtain based on the available model and forecast of the disturbance $w$, but makes no use of the fact that for a given $u$ and $w$, the system operator can often \emph{measure} $y$. Recently, motivated by applications in power systems, there has been a great effort to develop \emph{feedback} based optimization schemes that make use of the the measurement $y$ in real time in order to attempt to solve~\eqref{eq.opt}. 
%
\subsection{Real-time feedback optimization - online approximate gradient descent}
Let us introduce a soft-constrained version of problem~\eqref{eq.opt.orig}
\begin{align}\label{eq.opt.soft}
\begin{split}
    \min_{u\in\mc U} ~&~f(u) + g(y)\\
    \mathrm{subject~to}&~y = \pi(u,w),\\
    \end{split}
\end{align}
where $g$ is an appropriate convex, continuously differentiable strongly smooth penalty function {for the constraint $y \in \mathcal{Y}$}. The approximation~\eqref{eq.opt.soft} is justified in many applications, as the set $\mc Y$ often represents {desired} engineering constraints (i.e., voltage or line current limits in a power networks) for which small and/or infrequent violations carry no big consequences. The hard constraints on the inputs $\mc U$ are not relaxed, and will be enforced at all times. Next, we present the standard (measurement-based) gradient descent for~\eqref{eq.opt.soft}.
\vspace{-0.3cm}
\begin{center}
\begin{algorithm}[ht]
\SetAlgoLined
\KwIn{$k = 1,\, \tau>0,\,  u_1\in\mc U$ }
 \Iterate
 {\textbf{measure} $\displaystyle{ y_k = \pi(u_k,w)}$ \\[0.2em]
 $\displaystyle {d_{k} = \nabla f(u_k) + \partial\pi(u_k,w)\tr\nabla g(y_k)}$ \\[0.2em]
 $\displaystyle {u_{k+1} = \mathrm{Proj}_\mc U \left(u_{k}  - \tau d_k\right)} $ \\[0.2em]
 \textbf{apply} $u_k$ to the system\\[0.2em]
 $\displaystyle  k\leftarrow k+1$\\[0.1em]}
 \caption{Gradient Descent (GD)}\label{alg.primal.grad.exact}
\end{algorithm} 
\end{center}
\vspace{-0.8cm}
In Algorithm~\ref{alg.primal.grad.exact}, with a slight abuse of notation, we denote by\looseness=-1
\[
\partial \pi(\bar u,\bar w) := \left.\frac{\partial \pi(u,w)}{\partial u}\right|_{u = \bar u, w= \bar w}
\]
the Jacobian of $\pi$ with respect to $u$ evaluated at a point $\bar u,\bar w$. Under suitable assumptions on the step size $\tau$, Algorithm~\ref{alg.primal.grad.exact} converges to a KKT point of~\eqref{eq.opt.soft}. Note that in order to implement Algorithm~\ref{alg.primal.grad.exact}, we no longer need a model of the full map $\pi$, as $\pi(u_k,w)$ can be measured in real-time through $y$. Instead, we need only its Jacobian with respect to the 
decision variable $u$ (which still depends $\pi$ and $w$). 
In~\cite{Dallanese2016optimal,hauswirth2016projected}, the authors show promising results for various optimal power flow problems where algorithms similar to Algorithm~\ref{alg.primal.grad.exact} reach near-optimal solutions, reject time-varying disturbances despite not using the exact Jacobian; these robustness properties are precisely the well-known advantages of feedback over feedforward control. In this work, we study robustness for the simplest approximation of the Jacobian i.e.,
\begin{align}\label{eq.approx.lin}
\partial\pi(u,w) \approx \Pi ,\quad \forall\,u\in\mc U,  w\in\mc W.
\end{align}
Using~\eqref{eq.approx.lin}, Algorithm~\ref{alg.primal.grad.exact} becomes the Online Approximate Gradient (OAG) algorithm (Algorithm~\ref{alg.primal.grad}).
\begin{center}
\begin{algorithm}[ht]
\SetAlgoLined
\KwIn{ $k = 1,\, \tau>0,\,  u_1\in\mc U$ }
 \Iterate
 {\textbf{measure} $\displaystyle{ y_k \leftarrow \pi(u_k,w)}$ \\[0.2em]
 $\displaystyle {u_{k+1} = \mathrm{Proj}_\mc U \left(u_{k}  - \tau(\nabla f(u_k) + \Pi\tr\nabla g(y_k))\right)} $ \\[0.2em]
  \textbf{apply} $u_k$ to the system\\[0.2em]
 $\displaystyle  k\leftarrow k+1$\\[0.1em]}
 \caption{Online Approximate Gradient (OAG)}\label{alg.primal.grad}
\end{algorithm} 
\end{center}
\vspace{-0.8cm}

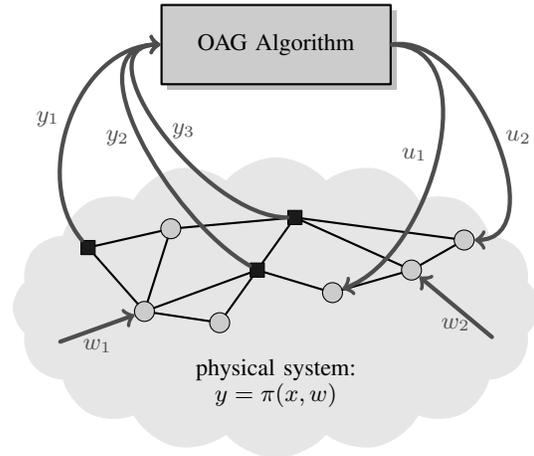
\begin{figure}[htbp]
\begin{center}
\begin{tikzpicture}[auto, thick]
  \node[cloud, fill=gray!20, cloud puffs=16, cloud puff arc= 100,
        minimum width=7cm, minimum height=4cm, aspect=1] at (0,-0.5) {};

  \foreach \place/\x in { {(-1.75,-0.55)/2},{(-1.4,0.55)/3},
    {(-0.75,-0.7)/4}, {(0.75,-0.3)/7}, 
    {(1.8,0)/8},{(2.5,0.4)/9}}
  \node[cgray] (a\x) at \place {};

  \foreach \place/\x in {{(-2.5,0.3)/1}, 
     {(-0.25,0)/5}, {(0.25,0.7)/6}}
  \node[qblack] (s\x) at \place {};

  \node[block_m] at (0,3) (controller) {OAG Algorithm};

  \path[thin] (s1) edge (a2);
  \path[thin] (s1) edge (a3);
  \path[thin] (a2) edge (a3);
  \path[thin] (a3) edge (s6);
  \path[thin] (a2) edge (a4);
  \path[thin] (s5) edge (s6);
  \path[thin] (s5) edge (a4);
  \path[thin] (s5) edge (a2);
  \path[thin] (s5) edge (a7);
  \path[thin] (s6) edge (a9);
  \path[thin] (s6) edge (a8);
  \path[thin] (a8) edge (a9);
  \path[thin] (a7) edge (a8);
 
   \path[gpath] (s1) to [out=130,in=180] node [left] {$y_1$} (controller.west);
   \path[gpath] (s5) to [out=150,in=180] node [left] {$y_2$} (controller.west);
   \path[gpath] (s6) to [out=180,in=180] node [right] {$y_3$} (controller.west);

   \path[gpath] (controller.east) to [out=0,in=20] node [left] {$u_1$} (a7);
   \path[gpath] (controller.east) to [out=10,in=0] node [right] {$u_2$} (a9);
   
   \node (w1) at (-3,-1) {};
   \node (w2) at (3,-1) {};

   \path[gpath] (w1) to  node [below] {$w_1$} (a2);
   \path[gpath] (w2) to  node [below] {$w_2$} (a8);
   \node (N) at (0,-1.5) {$\begin{array}{c}\text{physical system:} \\ y = \pi(x,w) \end{array}$};

\end{tikzpicture}
\caption{Online implementation of Algorithm~\ref{alg.primal.grad}: The update of $x$ is computed based on the measurement of ${y} = \pi(x,w)$ and applied to the system.}
\vspace{-0.6cm}
\label{fig.online.extragrad}
\end{center}
\end{figure}

Note that, Algorithm~\ref{alg.primal.grad} can be implemented in a ``online" fashion using the approximate Jacobian $\Pi$ and the system measurements $y$ (i.e., using feedback). No information on the real system model $\pi$ or the disturbance $w$ is required. This is illustrated in Figure~\ref{fig.online.extragrad}. In fact Algorithm~\ref{alg.primal.grad} can be interpreted as a hybrid between an optimization algorithm and a feedback controller that tries to \emph{steer} the system close to the optimal solution of Problem~\eqref{eq.opt.soft}. 
The simple approximation~\eqref{eq.approx.lin} is consistent with the fact that, as discussed in Section~\ref{sec.feedforward}, many system operators use linear models of their physical systems for large-scale optimization problems.

In the remainder of the work we will make some assumptions on the Jacobian $\partial\pi (u,w)$ in the set of interest $\mc U\times \mc W$ and we will try to analyze the stability and robustness properties of the OAG algorithm. In particular we are interested in the following questions
\begin{itemize}
\item Can we characterize the set of points, if any, to which the OAG algorithm converges?
\item Can we guarantee that the OAG algorithm is robustly stable (i.e., it converges for a large class of maps $\pi$)? 
\end{itemize}
\subsection{Characterizing the closed-loop equilibria of OAG}
Observe that, while our goal is to solve the non-convex optimization \eqref{eq.opt.orig} using feedback, the OAG algorithm  uses only an approximation $\Pi$ of $\partial \pi(u,w)$; we are therefore unlikely to converge to a true optimizer of \eqref{eq.opt.orig}, but instead will be able to converge to a point which is simultaneously consistent with both the OAG algorithm and the physical system. We formalize this idea as that of an \emph{online approximate solution}.


{\begin{definition}[Online approximate solution]\label{def.approximate}
Given $\Pi\in\R^{m\times n}$ and a disturbance $w\in\R^q$, a vector $\bar u = \bar u(w) \in\R^n$ is an online approximate solution of~\eqref{eq.opt.soft} if\looseness=-1
\begin{subequations}\label{Eq:approximateconditions}
\begin{align}
&  y = \pi(\bar u,w)\\
&\bar u \in\mc U, \label{feas1}\\
&-\nabla f(\bar u)- \Pi\tr\nabla g(y) \in\mc N_{\mc U}^{\bar u},
\label{improve}
\end{align}
\end{subequations}
where $\mc N_{\mc U}^{\bar u}$ is the normal cone of $\mc U$ at the point $\bar u$.
\end{definition}
In other words an online approximate solution is a {feasible} solution which \emph{would} be a KKT point of~\eqref{eq.opt.soft} if the linear model was locally accurate (i.e., $\partial \pi(\bar u,w) = \Pi$). This means that if the decision maker \emph{believes} her linearized model, she has no incentive to change the solution\footnote{This is assuming that the decision maker is after a KKT point, which is reasonable given that~\eqref{eq.opt} is a non-convex problem}. The quality of an online approximate solution depends entirely on the quality of the linear model.
}
Next, we provide a condition for which the OAG algorithm converges to an online approximate solution. To do so, we assume, without loss of generality, a specific structure for the set $\mc U$. 

\begin{assumption}\label{ass.U} The set $\mc U$ is partitioned as
$
\mc U = \mc U_1 \times \mc U_2 \times \mc U_3,
$
where $\mc U_1 := \R^{n_1}$, $\mc U_2 \subset \R^{n_2}$ is a box constraint i.e., $\mc U_2 := \{u_2\in\R^{n_2}\,|\, \underline u_{2,i} \le u_{2,i} \le \bar u_{2,i}, i=1,\dots,n_2 \}$ $-\infty\le\underline u_{2,i}  \le \overline u_{2,i}\le \infty$ and $\mc U_3 \subset \R^{n_3}$ is a general closed convex set and $n_1+n_2+n_3 = n$. 
\end{assumption}
Next, we define the set
\[
\mc P := \left\{ P \succ 0 \,\Bigg |\, P = \begin{bmatrix}P_1 & & \\ & P_2 & \\ && I_{n_3}\end{bmatrix}, \begin{array}{c}P_1\in \R^{n_1\times n_1} \\
P_2\in  \mathbb D_{n_2} 
\end{array}  \right\}
\]
where $\mathbb D_{n}$ is the set of $n\times n$ diagonal matrices.
\begin{proposition}[Convergence of OAG]\label{propo.convergence}
Let $\mc U$ satisfy Assumption~\ref{ass.U}. For $P\in\mc P$, assume that $F_w(u) := \nabla f(u) + \Pi\tr\nabla g(\pi(u,w))$ is $\rho-$strongly monotone and $L-$Lipschitz continuous w.r.t $\langle \cdot,\cdot\rangle_{P}$.  Then if $\tau < \frac{2\rho}{L^2}$, the OAG algorithm converges geometrically to the unique online approximate solution.
\end{proposition}
\begin{proof}
By the definition of normal cone, an online approximate solution $\bar u$ satisfies
$
F_w(\bar u)\tr (u-\bar u) \ge 0,\, \forall\,u\in\mc U\,.
$
or equivalently
$
\bar u = \mathrm{Proj}_{\mc U} (\bar u - \tau F_w(\bar u)),
$
and therefore solves $\vi{\mc U}{F_w}$ with respect to the standard Euclidean norm. By applying the definition of projection, it is easy to see that, under Assumption 1, if  $P\in\mc P$, then $\mathrm{Proj}_{\mc U}^P (u) = \mathrm{Proj}_{\mc U} (u)$ for all $u\in\R^n$. By Proposition~\ref{propo.conv.grad}, $\vi{\mc U}{F_w}$ has a unique solution $\bar u$ (an online approximate solution) and the OAG converges geometrically to $\bar u$ for any initial condition $u_1$. 
\end{proof}

\smallskip

With the following proposition, we show that the intuition behind an online approximate solution is indeed correct and we can bound the distance of an online approximate solution to a KKT point $u^\star$ of \eqref{eq.opt.soft} based on the Jacobian approximation error $\|\Pi - \partial\pi(u^\star,w))\|$.

\begin{proposition}[Approximation error]
If $P\in\mc P$ and $F_w(u) := \nabla f(u) + \Pi\tr\nabla g(\pi(u,w))$ is $\rho-$strongly monotone w.r.t $\langle \cdot,\cdot\rangle_{P}$, then
\[
\|\bar u - u^\star\|_P\le \frac{1}{\rho}\|(\Pi - \partial\pi(u^\star,w))\tr\nabla g(\pi(u^\star,w))\|_P
\]
where $\bar u$ is the unique online approximate solution satisfying \eqref{Eq:approximateconditions} and $u^\star$ a KKT point of \eqref{eq.opt.soft}.
\end{proposition}
\begin{proof}Let us define 
$
F_w^\star(u) := \nabla f(u) + \partial\pi(u,w)\tr\nabla g(\pi(u,w)).
$
The result follows from noting that $u^\star$ solves $\vi{\mc U}{F_w^\star}$ and applying~\cite[Theorem 1.14]{nagurney2013network}. 
\end{proof}



\section{Robust monotonicity of uncertain operators}

{ We are now interested in developing conditions to check whether \textemdash{} given a suitable choice of $\Pi$ \textemdash{} the OAG algorithm (Algorithm \ref{alg.primal.grad}) is robustly stable when implemented on the uncertain physical system described by $\pi(u,w)$.} To do so, we begin by abstracting the OAG algorithm by writing the $u$ update compactly as
\[
u_{k+1} = \mathrm{Proj}_\mc U \left(u_{k}  - \tau F_w(u_k) \right)
\]
where 
\[
F_w(u):= \nabla f(u) + \Pi\tr\nabla g(\pi(u,w)).
\] 
According to Proposition~\ref{propo.conv.grad}, if $F_w$ is strongly monotone and Lipschitz continuous, then the OAG algorithm converges geometrically for $\tau < 2\rho/L^2$. { Instead of directly modelling $\pi(u,w)$ as uncertain, we will consider uncertainty on the map $F_w$, and develop conditions under which (strong) monotonicity can be guaranteed robustly with respect to this uncertainty. Inspired by Proposition \ref{propo.jen.jac}, we will parametrize }uncertainty on $F_w$ by defining an uncertainty set $\mc J$, and imposing that $\partial _CF_w(u)\subset \mc J$, for all $u\in\mc U$ and all $w\in\mc W$. We investigate the simple case in which $\mc J$ is a polytope and the more interesting case in which $\mc J$ is parametrized by a Linear Fractional Transformation (LFT).  

\subsection{Polytopic uncertainty in the Jacobian}

We begin with the simple case where 
\begin{equation}\label{Eq:JPoly}
 \mc J^{\text{poly}} := \mathrm{co}\left\{J_i,~i = 1,\dots,\nu\right \}.
\end{equation}
Then, given the set $\mc U$, we define the set of functions
\begin{multline}\label{Eq:FPoly}
\mc F^{\text{poly}} := \{F_w \,|\, \partial_C F_w(u)\subseteq   \mc J^{\text{poly}},  \forall \,u\in\mc U  \}.
\end{multline}
The following proposition provides a numerical test to guarantee that strong monotonicity holds robustly for all $F_w\in\mc F^{\text{poly}}$.
\begin{proposition}\label{propo.polytope}
Given $P\succ 0$ and a constant $\rho > 0$, the following two statements are equivalent:
\begin{enumerate}[(i)]
\item all operators $F_w\in\mc F^{\text{poly}}$, are $\rho-$strongly monotone w.r.t $\langle \cdot,\cdot\rangle_{P}$ on the set $\mc U$;
\item the following Matrix Inequality holds true
\begin{align}\label{eq.lmi.polytope}
\frac{1}{2}\left[J_i\tr P + P J_i\right] \succ \rho P,\quad i=1,\dots,\nu.
\end{align}
where $J_i$ as in~\eqref{Eq:JPoly}.
\end{enumerate}
Moreover, if $\rho = 0$, then the the preceding statements are equivalent with ``monotone'' replacing ``strongly monotone''.

\end{proposition}
The proof is simple and is omitted for reasons of space. Since~\eqref{eq.lmi.polytope} is a LMI in $P$, Proposition~\ref{propo.polytope} allows to test the hypothesis of Proposition~\ref{propo.convergence} and thus can be used to guarantee that the OAG algorithm converges robustly for all operators $F_w \in \mathcal{F}^{\rm poly}$. In Section~\ref{sec.application.optimization} we will show how Proposition~\ref{propo.polytope} can be applied to the operator that arises from the OAG algorithm with a prototypical soft-constrained optimization problem of the form~\eqref{eq.opt.soft}. Not surprisingly, the number of constraints in the LMI~\eqref{eq.lmi.polytope} can be very large in practical problems. For this reason in the next section we will consider an arguably better way to parametrize uncertainty in the Jacobian, which, at the expense of an increased modeling effort, leads to a more elegant test involving a single LMI.

\subsection{LFT uncertainty in the Jacobian}\label{sec.lft}

We now consider a different parametrization for the uncertainty set of the Clarke generalized Jacobian of $F_{w}$. Given a set of matrices $\mathbf \Delta\subset \R^{s\times z}$ and fixed matrices $A\in\R^{n\times n}, B\in\R^{n\times s}, C\in\R^{z\times n}$ and $D\in\R^{z\times s}$, define
\begin{align*}
\mc J^{\text{lft}} := 
\left\{A + B\Delta(I_z-D\Delta)^{-1} C\,:\, \Delta \in \mathbf \Delta \right \}\,,
\end{align*}
where we assume that $(I_z-D\Delta)$ is invertible for all $\Delta \in\mathbf \Delta$. In addition, we suppose we have access to a convex cone of matrices $\mathbf\Theta\subset\R^{(s+z)\times (s+z)}$ such that
\begin{align}\label{eq.iqc}
p = \Delta\, q,~ \Delta \in \mathbf \Delta \implies 
\begin{bmatrix}
q\\
p
\end{bmatrix}\tr \Theta \begin{bmatrix}
q\\
p
\end{bmatrix} \ge 0,  \,\forall\,
\Theta\in\mathbf\Theta \,.
\end{align}
%
%
%
The parametrization of uncertainty and the positivity criteria \eqref{eq.iqc} is borrowed from the literature on robust control \cite{GED-FP:00,AM-AR:97,CS-SW:15,JV-CWS-HK:16}. This parametrization might seem unnatural at first but it is extremely powerful in modeling a large set of common uncertainty classes; unfortunately it requires a steep learning curve to get accustomed to. 
As we did in \eqref{Eq:FPoly}, we let $\mathcal{F}^{\rm lft}$ denote the set of functions such that $\partial_C F_w(u) \subset \mc J^{\text{lft}}$ for all $u \in \mathcal{U}$. The next result shows that robust monotonicity in this set can be tested with a single LMI.
 
\begin{proposition}\label{propo.lft}
If there exists $P \succ 0$ and $\Theta\in\mathbf\Theta$ such that
\begin{align}\label{eq.test}
\begin{split}
\begin{bmatrix}
A_\rho \tr P + P A_\rho & P B\\
B\tr P & \vzeros
\end{bmatrix}
-
\begin{bmatrix}
C & D\\
\vzeros &  I_s 
\end{bmatrix} \tr
\Theta
\begin{bmatrix}
C & D\\
\vzeros &  I_s 
\end{bmatrix} 
\succcurlyeq 0\\
\end{split}
\end{align}
where $A _\rho = A -\rho I_n$, then all operators $F_w \in \mathcal{F}^{\rm lft}$ are $\rho-$strongly monotone w.r.t $\langle \cdot,\cdot\rangle_{P}$ over $\mc U$, or are simply monotone if~\eqref{eq.test} holds for $\rho = 0$.
\end{proposition}
\begin{proof}
Given any $x\in \R^n$ and $\Delta \in \boldsymbol{\Delta}$, let us define $q = Cx+Dp$ and $p = \Delta q$. Pre and post multiplying~\eqref{eq.test} by $[x\tr p\tr ]\tr$, we obtain 
\begin{align}\label{eq.inequality.pq}
\begin{split}
\begin{bmatrix}
x \\ p
\end{bmatrix} \tr
\begin{bmatrix}
A_\rho \tr P + P A_\rho & P B\\
B\tr P & \vzeros
\end{bmatrix}
\begin{bmatrix}
x \\ p
\end{bmatrix}  
-
\begin{bmatrix}
q \\ p
\end{bmatrix}  \tr
\!\!\Theta
\begin{bmatrix}
q \\ p
\end{bmatrix}  
\ge 0
\end{split}
\end{align}
Since $p = \Delta q$,~\eqref{eq.iqc} implies that the first term in~\eqref{eq.inequality.pq} must be nonnegative. Therefore,
\begin{align}
\begin{split}
x\tr P A_\rho  x  & + x\tr P B p  \ge 0,\\
\forall x\in \R^n \text{ and }~ &q = Cx + Du\,,\,\, p = \Delta q.
\end{split}
\end{align}
We conclude that $p  = \Delta(I_{z}-D\Delta)^{-1}Cx$. Substituting into the first inequality, we obtain
\begin{equation}\label{Eq:JDelta}
\frac{1}{2}\left[J_{\Delta}\tr P + PJ_{\Delta} \right] \succeq \rho P 
\end{equation}
where $J_{\Delta} = A + B\Delta(I_{z}-D\Delta)^{-1} C \in \mathcal{J}^{\rm lft}$. Hence, all elements of $\mathcal{J}^{\rm lft}$ satisfy \eqref{Eq:JDelta}, and therefore all functions $F_w \in \mathcal{F}^{\rm lft}$ are $\rho$-strongly monotone w.r.t $\langle \cdot,\cdot\rangle_{P}$.
\end{proof}
Since~\eqref{eq.test} is a LMI in $P$ and $\Theta$, Proposition~\ref{propo.lft} allows to verify the Proposition~\ref{propo.convergence} numerically. Thus it can be used to guarantee that the OAG algorithm converges robustly for all operators $F_w \in \mathcal{F}^{\rm lft}$
\begin{remark}
Proposition~\ref{propo.convergence} requires Lipschitz continuity of $F_w$ w.r.t $\langle \cdot,\cdot\rangle_{P}$ to guarantee stability. Note that, in our case $F_w$ is robustly $L$-Lipschitz w.r.t $\langle \cdot,\cdot\rangle_{P}$ if, for all $J\in\mc J$, $J\tr P J - L^2 P \preccurlyeq 0$. This can be tested using LMIs both for $\mc J^{\text{poly}}$ and $\mc J^{\text{lft}}$. This analysis is omitted for reasons of space.
\end{remark}
\subsection{Quick recipes for LFT modeling}\label{Sec.cookbook}
The hard work in applying Proposition \ref{propo.lft} typically comes in writing down a useful cone of matrices $\boldsymbol{\Theta}$ such that the positivity condition in \eqref{eq.iqc} holds. Luckily, there are standard recipes for doing this for some practically important uncertainty sets $\boldsymbol{\Delta}$. While LFT modeling is applicable to uncertain nonlinear operators (see \cite[Chapter 6]{CS-SW:15} for a complete treatment), in this section we will limit the treatment to uncertain matrices as they suffice to model the uncertainty in $\partial_{C} F_w$ described in Section~\ref{sec.lft}.
\subsubsection{Unstructured, norm-bounded uncertainty}
Given $\gamma \geq 0$, let
$
\boldsymbol{\Delta}_{\rm u,nb}(\gamma) :=  \setdef{\Delta \in \R^{s \times z}}{\|\Delta\|_2 = \sigma_{\rm max}(\Delta) \leq \gamma}
$
which is the set of unstructured matrices with induced norm less than or equal to $\gamma$. A cone $\boldsymbol{\Theta}$ that achieves the required positivity condition in \eqref{eq.iqc} is 
\[
\begin{aligned}
\boldsymbol{\Theta}_{\rm u,nb} &= \Big\{\theta \left[\begin{smallmatrix}I_s & 0\\ 0 & -\frac{1}{\gamma^2} I_z\end{smallmatrix}\right]\,\,\Big|\,\,\theta \geq 0\Big\}\,.
\end{aligned}
\]
To see this, note that if $\Delta \in \boldsymbol{\Delta}_{\rm u,nb}$, then $\|p\|_2^2 = \|\Delta q\|_2^2 \leq \gamma^2 \|q\|^2$, and therefore $\theta(\|q\|_2^2 - \frac{1}{\gamma^2} \|p\|_2^2) \geq 0$ for any $\theta \geq 0$.

\subsubsection{Repeated scalar norm-bounded uncertainty} Given $\gamma \geq 0$ let $\boldsymbol{\Delta}_{\rm rs,nb}(\gamma) :=  \setdef{\Delta = \delta I}{\delta \in \R,\,\, |\delta| \leq \gamma}$ denote the set of uniform diagonal matrices with diagonal entries bounded in magnitude by $\gamma$. A cone $\boldsymbol{\Theta}$ that works is
\[
\begin{aligned}
\boldsymbol{\Theta}_{\rm rs,nb} &= \Big\{\Theta = \left[\begin{smallmatrix}\Phi & \Psi\\ \Psi \tr & -\frac{1}{\gamma^2} \Phi \end{smallmatrix}\right]\,\,\Big|\,\,
\Phi \succcurlyeq 0,\,\Psi = -\Psi \tr\Big\}\,.
\end{aligned}
\]
as may be verified by direct computation. Note that since we know more about the structure of the uncertainty, we can use a larger cone of matrices $\boldsymbol{\Theta}$; this reduces conservatism.

\subsubsection{Unstructured monotone and Lipschitz uncertainty} Given $\rho, L \in \R$ satisfying $0 \leq \rho \leq L < \infty$, let $\boldsymbol{\Delta}_{\mathrm{u},\rho L}$ denote the set of matrices ${\Delta}\in\R^{s\times s}$ such that $\rho I\preccurlyeq \Delta \preccurlyeq L I $. A cone that works for this case is
\[
\begin{aligned}
\boldsymbol{\Theta}_{\mathrm{u},\rho L} &= \Big\{\Theta = \varphi\left[\begin{smallmatrix}-2\rho L & \rho+L\\ \rho+L & -2 \end{smallmatrix}\right] \otimes I_s\,\,\Big|\,\,
\varphi \geq 0\Big\}\,.
\end{aligned}
\]
as may be verified again by direct calculation.

\subsubsection{Repeated scalar monotone and Lipschitz uncertainty} Given $\rho, L \in \R$ satisfying $0 \leq \rho \leq L < \infty$, let $\boldsymbol{\Delta}_{\mathrm{rs},\rho L}$ denote the set of diagonal matrices $\Delta = \delta I$ with $\rho \le\delta\le L$. A cone that works for this case is
\[
\begin{aligned}
\boldsymbol{\Theta}_{\mathrm{rs},\rho L} &= \Big\{\Theta = \varphi\left[\begin{smallmatrix}-2\rho L\Phi & (\rho+L)\Phi\\ (\rho+L)\Phi & -2\Phi \end{smallmatrix}\right]\,\,\Big|\,\, \Phi\succcurlyeq 0 
 \Big \}\,.
\end{aligned}
\]

\subsubsection{Block-structured uncertainty} Consider now the block-diagonal uncertainty set
\begin{align*}
\boldsymbol{\Delta}_{\rm blk} &:= \big\{\Delta = \mathrm{blkdiag}(\Delta_1,\ldots,\Delta_{r})\,\,\big|\,\,
  \Delta_{i} \in \R^{s_{i} \times z_{i}}, \\
 & \Delta_{i} \in \boldsymbol{\Delta}_{\rm u,nb} \,\,\text{or}\,\, \boldsymbol{\Delta}_{\rm rs,nb} \,\,\text{or}\,\, \boldsymbol{\Delta}_{\mathrm{u},\rho L}\,\, \text{or}\,\, \boldsymbol{\Delta}_{\mathrm{rs},\rho L}\big\}\,.
\end{align*}
where each block satisfies one of the previous criteria. Then the previous cones may be used individually for each corresponding block of the uncertainty.


%

\section{Application to Feedback Optimization}\label{sec.application.optimization}

Consider the following optimization problem
\begin{align}\label{eq.opt1}
\min_{u\in\mc U} ~& u\tr H u + h\tr u + \eta\sum_{i=1}^m \left\{ \max(0, \underline {y_i}-y_i, y_i-\overline {y_i})\right\}^2\nonumber \\
\mathrm{s.t.} ~& y = \pi(u,w).
\end{align}
in which the $\max$ functions encode soft versions of the constraints $y_i\in \left[\underline {y_i},\overline {y_i}\right]$ and $H\succ0$.
Given an approximator $\Pi$ of $\partial \pi$, if we run the OAG algorithm applied to~\eqref{eq.opt1}, the approximate gradient $F_w$ takes the form
\begin{align}\label{eq.Fw.opt}
F_w(u) = H u + h + \eta\,\Pi\tr s_{\underline y,\overline y}(\pi(u,w))\,,
\end{align}
where $s_{\underline{y},\overline{y}}$ is the (vectorized) soft-thresholding function (linear with unit slope for $y_i\not\in \left[\underline {y_i},\overline {y_i}\right]$, zero otherwise).
%
%
%
%
We now illustrate how to use the LMI conditions of Propositions~\ref{propo.polytope} and~\ref{propo.lft}  to guarantee robust stability of the OAG algorithm applied to~\eqref{eq.opt1}.
\subsection{Polytopic uncertainty}
We begin with polytopic uncertainty to show how this parametrization of the uncertainty, despite being the most intuitive, can quickly lead to an intractable number of constraints. Suppose the Jacobian $\partial \pi$ of $\pi(\cdot,w)$ lies in the convex hull of a set of known matrices
\[
\partial \pi(u,w)\in\mathrm{co}\left\{ \tilde \Pi_{i}, \,i=1,...,\nu  \right\}, \quad \forall\, u\in\mc U, w\in\mc W,
\]
then, for all $w\in \mc W$  the Clarke generalized Jacobian of the approximate gradient $F_w$ defined in~\eqref{eq.Fw.opt} is satisfies $\partial_CF_w \in \mc J^{\text{poly}}$, with 
\[
\mc J^{\text{poly}} := \mathrm{co}\left\{H+\eta\,\Pi\tr Q_j\tilde\Pi_i,~ i=1,...,\nu,  j=1,...,2^m\right \},
\]
where the matrices $Q_j$ are diagonal with all possible combinations of diagonal elements in $\{0,1\}$. Let $M_{ji} = Q_j\tilde{\Pi}_{i}$, from Propositions~\ref{propo.convergence} and~\ref{propo.polytope} we know that, under Assumption~\ref{ass.U}, if we find $P\in\mc P$ such that 
\begin{align}\label{eq.mono.LMI.2}
 (H +\eta\, \Pi \tr M_{ji}) P  + P( H +\eta\, \Pi\tr M_{ji})\tr \succcurlyeq \rho P
\end{align}
for $i=1,...,\nu$, $j=1,...,2^m$, then the OAG algorithm converges to the unique online approximate solution. 
Even for this simple example, we obtain condition \eqref{eq.mono.LMI.2} with $p\cdot2^m$ LMI constraints, which become intractable even for fairly low-dimensional problems. 


\subsection{LFT uncertainty}\label{sec.application.optimization.lft}
Let us now consider a different parametrization for the uncertainty in the map $\pi$. Suppose that for all $w\in\mc W$, $u\in\mc U$ 
\begin{align}\label{eq.jacobian}
\partial \pi(u,w) = \Pi_\text{nom} + \Delta_\pi(u,w),~~\|\Delta_\pi(u,w)\|\le \gamma\,,
\end{align}
where $\Pi_{\rm nom}$ is a nominal value or ``best guess'' for $\partial \pi$ at normal operating conditions. By differentiating \eqref{eq.Fw.opt}, one may deduce that for all $w\in \mc W$  the Clarke generalized Jacobian of the mapping $F_w$ defined in~\eqref{eq.Fw.opt} satisfies $\partial_CF_w \in \mc J^{\text{lft}}$, where 
\begin{align*}
\mc J^{\text{lft}} &:= \{H  + \eta \Pi\tr\Delta_q ( \Pi_\text{nom} + \Delta_\pi)\,|\,\\
 &\qquad \Delta_q\in\mathbb D_m,\, 0 \preccurlyeq \Delta_q\preccurlyeq I_m,\, \|\Delta_\pi\|\le \gamma   \}.
\end{align*}
One can verify that each element $J$ of $\mc J^{\text{lft}}$ can be written as
$
J  = A + B(I-\Delta D)^{-1}\Delta C
$
with
\begin{align}\label{eq.abcd}
\left[\begin{array}{c|c}
A & B\\\hline
C & D
\end{array}\right] = \left[\begin{array}{c|cc}
H & \eta\, \Pi\tr & \vzeros \\ \hline
\Pi_\text{nom}  &\vzeros &  I_m \\
I_n & \vzeros & \vzeros 
\end{array}\right]
\end{align}
and 
$
\Delta =\mathrm{blkdiag} ( \Delta_1,  \Delta_2) 
$
with $\Delta _1 =  \Delta_q$ and $\Delta_2 = \Delta_\pi$. 
Using the recipes from Section~\ref{Sec.cookbook}, since the $j^\text{th}$ diagonal element of $\Delta_1$  is an independent scalar uncertainty in $[0,1]$ and $\Delta_2$ is unstructured and norm-bounded, we can show that all $\Delta\in\mathbf \Delta$ satisfy~$\eqref{eq.iqc}$ with $\mathbf \Theta$ of the form
 \begin{align}\label{eq.theta}
 \mathbf\Theta:=
\left\{
\sum_{j=1}^m
\Theta_j \,\Bigg|\, \varphi_j \geq 0, ~ \theta\geq 0
\right\},
\end{align}
with 
\[
\Theta_j := \left[\begin{smallmatrix}
\vzeros & \vzeros &\varphi_j\, \mathrm  e_j\mathrm e_j\tr & \vzeros \\
\vzeros &\frac{\theta}{m} I_n& \vzeros& \vzeros \\
\varphi_j\, \mathrm  e_j\mathrm e_j\tr & \vzeros & -2 \varphi_j\, \mathrm  e_j\mathrm e_j\tr & \vzeros\\
\vzeros & \vzeros & \vzeros& - \frac{\theta}{m\gamma^2} I_m \\
\end{smallmatrix}\right],
\]
where $\mathrm e_j$ is the $j^\text{th}$ canonical vector in $\R^m$. 
From Propositions~\ref{propo.convergence} and~\ref{propo.polytope} we know that, under Assumption~\ref{ass.U}, if we find $P\in\mc P$ and $\Theta\in\mathbf \Theta$ such that the single LMI~\eqref{eq.test} is satisfied, then the OAG algorithm converges to the unique online approximate solution.

\section{Examples and applications}\label{sec.examples}

In this section we illustrate the robustness properties of the OAG algorithm on a numerical example and on a timely engineering application in the context of power systems. 
\subsection{OAG numerical example}
\begin{figure}[htbp]
\begin{center}
\input{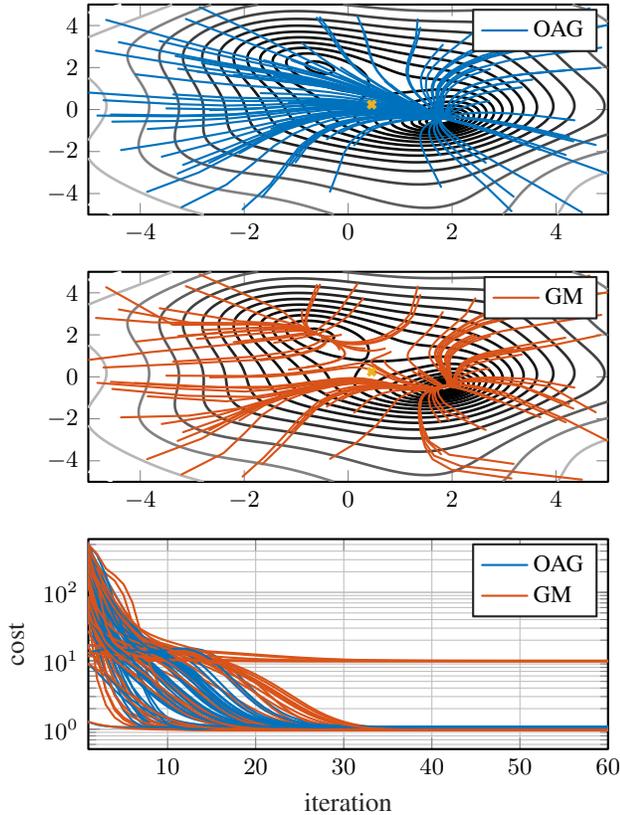}
\vspace{-0.3cm}
\caption{Comparison of the OAG algorithm (Algorithm~\ref{alg.primal.grad}) and the standard gradient method (GM Algorithm,~\ref{alg.primal.grad.exact}). The OAG algorithm converges to the unique online approximate solution, whose quality depends on the choice of $\Pi$. The yellow (\textcolor{mycolor2}{\textbf x}) marks the solution of the approximate convex program based on the linearization (i.e. the problem obtained by substituting $y = \Pi u$ in the constraints of~\eqref{eq.opt.p.example}). Clearly, this na\"ive solution is greatly outperformed by the OPE algorithm. There is no visible difference between the OPE algorithm and the GM in the number of iterations needed for convergence.}
\vspace{-0.6cm}
\label{fig.online.extragrad.example.primal1}
\end{center}
\end{figure}

In order to validate Propositions~\ref{propo.convergence} and~\ref{propo.polytope}, we consider the following physical system 
\begin{align*}
y = \pi(u,w) = \left\{ \begin{array}{l}
 u_1 + u_2\\
 w_1\,\sin(u_1) - u_1 +  w_2\,\cos(u_2) + u_2
\end{array}\right.
\end{align*}
where $w\in\mc W = [0,1]^2$ is an unknown disturbance. Given the set $\mc U = [-5,5]^2$. We are interested in the following optimization problem
\begin{align}\label{eq.opt.p.example}
\begin{split}
    \min_{u\in\mc U} ~&~u\tr Q_1 u + c_1\tr u + y\tr Q_2 y + c_2\tr y\\
    ~&~ y = \pi(u,w)
\end{split}
\end{align}
with $Q_1 = I, Q_2 = 10I, c_1 = [0,-9]\tr, c_2 = [-10,9]\tr$. We choose $\Pi =  \left[\begin{smallmatrix}
1 & 1 \\
-1 & 1
\end{smallmatrix}\right]$ to run OAG algorithm.
It is easy to see that $\partial F_w\in\mc J^{\text{poly}}$, where $J\in\mc J^{\text{poly}}=I_2+10\, \Pi\tr\tilde\Pi$ with
\[
\tilde \Pi \in \mathrm{co}\Bigg\{ 
\underbrace{\begin{bmatrix}
1 & 1 \\
0 & 0
\end{bmatrix}}_{\tilde\Pi_1},
\underbrace{\begin{bmatrix}
1 & 1 \\
0 & 2
\end{bmatrix}}_{\tilde\Pi_2},
\underbrace{\begin{bmatrix}
1 & 1 \\
-2 & 2
\end{bmatrix}}_{\tilde\Pi_3},
\underbrace{\begin{bmatrix}
1 & 1 \\
-2 & 0
\end{bmatrix}}_{\tilde\Pi_4}
\Bigg\}.
\]
%
%
By inspection, we observe that 
$
\Pi\tr  \tilde\Pi_i +  \tilde\Pi_i\tr  \Pi\succcurlyeq 0,~i = 1,...,4
$  
and therefore~\eqref{eq.lmi.polytope} is satisfied with $P=I$, $\rho =1$. By Proposition~\ref{propo.convergence}  the OAG algorithm converges to an online approximate solution of~\eqref{eq.opt.p.example} for any $w\in[0,1]^2$.
The evolution of the OAG algorithm (for $w = [1,1]\tr$), which uses the linear model and real-time feedback but does not have information on $w$ and $\pi$ is shown if Figure~\ref{fig.online.extragrad.example.primal1}, where it is compared with evolution the standard gradient method (Algorithm~\ref{alg.primal.grad.exact}) (that uses full information about $w$ and $\pi$) for one hundred random initial conditions. Both algorithms were implemented with a step-size $\tau=0.01$.
As we observe in Figure~\ref{fig.online.extragrad.example.primal1}, since the operator $F_w$ is strongly monotone, the OAG algorithm always converges to a single point while the gradient method can converge to different local minima. 

%

\subsection{Robust feedback optimization of a distribution feeder}
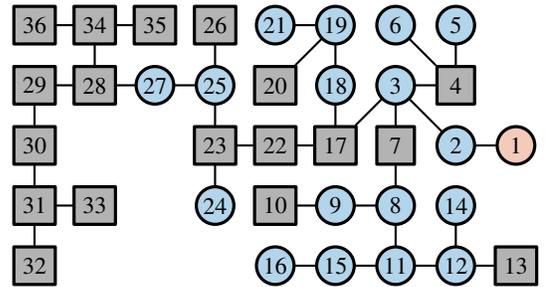
\begin{figure}[htb]
\begin{center}
\definecolor{mycolor1}{rgb}{0.00000,0.44700,0.74100}%
\definecolor{mycolor2}{rgb}{0.85000,0.32500,0.09800}%
\definecolor{mycolor3}{rgb}{0.92900,0.69400,0.12500}%

\newcounter{myval}

\begin{tikzpicture}[scale=0.8,vertex/.style={draw,circle,inner sep=0.1pt}, arc/.style={draw,thick,-}]


\foreach [count=\i] \coord in {
(0,0),
(-1,0),
(-2,1),
(-1,1),
(-1,2),
(-2,2),
(-2,0),
(-2,-1),
(-3,-1),
(-4,-1),
(-2,-2),
(-1,-2),
(0,-2),
(-1,-1),
(-3,-2),
(-4,-2),
(-3,0),
(-3,1),
(-3,2),
(-4,1),
(-4,2),
(-4,0),
(-5,0),
(-5,-1),
(-5,1),
(-5,2),
(-6,1),
(-7,1),
(-8,1),
(-8,0),
(-8,-1),
(-8,-2),
(-7,-1),
(-7,2),
(-6,2),
(-8,2)
}{
\setcounter{myval}{\i};
\def\mycol{mycolor1}
\IfSubStr{|1|} {|\arabic{myval}|}{\def\mycol{mycolor2}}{}
\IfSubStr{|4|7|10|13|17|20|22|23|26|28|29|30|31|32|33|34|35|36|} 
{|\arabic{myval}|}
{\def\mycol{black};
\node[draw, rectangle,minimum width= 0.5cm,minimum height= 0.5cm, very thick, fill=\mycol!30] (p\i) at \coord {\i};}
{\node[vertex, minimum width=0.5cm ,very thick, fill=\mycol!30] (p\i) at \coord {\i};}
}

\foreach [count=\r] \row in 
{
{0,1,0,0,0,0,0,0,0,0,0,0,0,0,0,0,0,0,0,0,0,0,0,0,0,0,0,0,0,0,0,0,0,0,0,0},
{0,0,1,0,0,0,0,0,0,0,0,0,0,0,0,0,0,0,0,0,0,0,0,0,0,0,0,0,0,0,0,0,0,0,0,0},
{0,0,0,1,0,0,1,0,0,0,0,0,0,0,0,0,1,0,0,0,0,0,0,0,0,0,0,0,0,0,0,0,0,0,0,0},
{0,0,0,0,1,1,0,0,0,0,0,0,0,0,0,0,0,0,0,0,0,0,0,0,0,0,0,0,0,0,0,0,0,0,0,0},
{0,0,0,0,0,0,0,0,0,0,0,0,0,0,0,0,0,0,0,0,0,0,0,0,0,0,0,0,0,0,0,0,0,0,0,0},
{0,0,0,0,0,0,0,0,0,0,0,0,0,0,0,0,0,0,0,0,0,0,0,0,0,0,0,0,0,0,0,0,0,0,0,0},
{0,0,0,0,0,0,0,1,0,0,0,0,0,0,0,0,0,0,0,0,0,0,0,0,0,0,0,0,0,0,0,0,0,0,0,0},
{0,0,0,0,0,0,0,0,1,0,1,0,0,0,0,0,0,0,0,0,0,0,0,0,0,0,0,0,0,0,0,0,0,0,0,0},
{0,0,0,0,0,0,0,0,0,1,0,0,0,0,0,0,0,0,0,0,0,0,0,0,0,0,0,0,0,0,0,0,0,0,0,0},
{0,0,0,0,0,0,0,0,0,0,0,0,0,0,0,0,0,0,0,0,0,0,0,0,0,0,0,0,0,0,0,0,0,0,0,0},
{0,0,0,0,0,0,0,0,0,0,0,1,0,0,1,0,0,0,0,0,0,0,0,0,0,0,0,0,0,0,0,0,0,0,0,0},
{0,0,0,0,0,0,0,0,0,0,0,0,1,1,0,0,0,0,0,0,0,0,0,0,0,0,0,0,0,0,0,0,0,0,0,0},
{0,0,0,0,0,0,0,0,0,0,0,0,0,0,0,0,0,0,0,0,0,0,0,0,0,0,0,0,0,0,0,0,0,0,0,0},
{0,0,0,0,0,0,0,0,0,0,0,0,0,0,0,0,0,0,0,0,0,0,0,0,0,0,0,0,0,0,0,0,0,0,0,0},
{0,0,0,0,0,0,0,0,0,0,0,0,0,0,0,1,0,0,0,0,0,0,0,0,0,0,0,0,0,0,0,0,0,0,0,0},
{0,0,0,0,0,0,0,0,0,0,0,0,0,0,0,0,0,0,0,0,0,0,0,0,0,0,0,0,0,0,0,0,0,0,0,0},
{0,0,0,0,0,0,0,0,0,0,0,0,0,0,0,0,0,1,0,0,0,1,0,0,0,0,0,0,0,0,0,0,0,0,0,0},
{0,0,0,0,0,0,0,0,0,0,0,0,0,0,0,0,0,0,1,0,0,0,0,0,0,0,0,0,0,0,0,0,0,0,0,0},
{0,0,0,0,0,0,0,0,0,0,0,0,0,0,0,0,0,0,0,1,1,0,0,0,0,0,0,0,0,0,0,0,0,0,0,0},
{0,0,0,0,0,0,0,0,0,0,0,0,0,0,0,0,0,0,0,0,0,0,0,0,0,0,0,0,0,0,0,0,0,0,0,0},
{0,0,0,0,0,0,0,0,0,0,0,0,0,0,0,0,0,0,0,0,0,0,0,0,0,0,0,0,0,0,0,0,0,0,0,0},
{0,0,0,0,0,0,0,0,0,0,0,0,0,0,0,0,0,0,0,0,0,0,1,0,0,0,0,0,0,0,0,0,0,0,0,0},
{0,0,0,0,0,0,0,0,0,0,0,0,0,0,0,0,0,0,0,0,0,0,0,1,1,0,0,0,0,0,0,0,0,0,0,0},
{0,0,0,0,0,0,0,0,0,0,0,0,0,0,0,0,0,0,0,0,0,0,0,0,0,0,0,0,0,0,0,0,0,0,0,0},
{0,0,0,0,0,0,0,0,0,0,0,0,0,0,0,0,0,0,0,0,0,0,0,0,0,1,1,0,0,0,0,0,0,0,0,0},
{0,0,0,0,0,0,0,0,0,0,0,0,0,0,0,0,0,0,0,0,0,0,0,0,0,0,0,0,0,0,0,0,0,0,0,0},
{0,0,0,0,0,0,0,0,0,0,0,0,0,0,0,0,0,0,0,0,0,0,0,0,0,0,0,1,0,0,0,0,0,0,0,0},
{0,0,0,0,0,0,0,0,0,0,0,0,0,0,0,0,0,0,0,0,0,0,0,0,0,0,0,0,1,0,0,0,0,1,0,0},
{0,0,0,0,0,0,0,0,0,0,0,0,0,0,0,0,0,0,0,0,0,0,0,0,0,0,0,0,0,1,0,0,0,0,0,0},
{0,0,0,0,0,0,0,0,0,0,0,0,0,0,0,0,0,0,0,0,0,0,0,0,0,0,0,0,0,0,1,0,0,0,0,0},
{0,0,0,0,0,0,0,0,0,0,0,0,0,0,0,0,0,0,0,0,0,0,0,0,0,0,0,0,0,0,0,1,1,0,0,0},
{0,0,0,0,0,0,0,0,0,0,0,0,0,0,0,0,0,0,0,0,0,0,0,0,0,0,0,0,0,0,0,0,0,0,0,0},
{0,0,0,0,0,0,0,0,0,0,0,0,0,0,0,0,0,0,0,0,0,0,0,0,0,0,0,0,0,0,0,0,0,0,0,0},
{0,0,0,0,0,0,0,0,0,0,0,0,0,0,0,0,0,0,0,0,0,0,0,0,0,0,0,0,0,0,0,0,0,0,1,1},
{0,0,0,0,0,0,0,0,0,0,0,0,0,0,0,0,0,0,0,0,0,0,0,0,0,0,0,0,0,0,0,0,0,0,0,0},
{0,0,0,0,0,0,0,0,0,0,0,0,0,0,0,0,0,0,0,0,0,0,0,0,0,0,0,0,0,0,0,0,0,0,0,0},
}
{
 \foreach [count=\c] \cell in \row
    {
    \ifnum\cell=1%
        \draw[arc] (p\r) edge (p\c);
    \fi
    }
}
\end{tikzpicture}
\caption{IEEE $37$-node feeder. Node 1 is the Point of Common Coupling (PCC). All other nodes are connected to a load and a voltage sensor. The square nodes are equipped with PV systems. The OAG algorithm is used to optimally decide curtailment of the PV systems in real-time in order to limit over-voltage.}
\vspace{-0.6cm}
\label{fig.test.feeder}
\end{center}
\end{figure}

In this section we illustrate how the LFT robustness test introduced in Section~\ref{sec.application.optimization.lft} can be used to certify robust stability of the OAG algorithm used to optimally manage the operation of a distribution feeder with high renewable penetration. The feeder, whose details can be found in~\cite{Dallanese2016optimal}, is illustrated in Figure~\ref{fig.test.feeder}. We simulate ten hours using real data from Anatolia, CA, USA, for solar irradiance and load consumption with granularity of one second. 
Let $u\in\R^{36}$ collect all controllable active and reactive power injections at the PV buses, $w\in\R^{70}$ collects all uncontrollable loads and power injections (active and reactive) at every node and $y\in\R^{35}$ collects the voltage magnitude at every node (except the PCC). The lower and upper voltage limits are $\underline u = 0.95$ p.u and $\overline u = 1.05$ p.u. Given the available active power for every PV $\{p^{\max}_i\}_{i=1}^{18}$ and the rated apparent power of every PV inverter $\{s^{\text {rated}}_i\}_{i=1}^{18}$, we define the closed convex power constraint sets 
\[
\mc U_i := \{u_i = [p_i,q_i]\tr\,| \, 0 \le p_i\le p^{\max}_i,\, q_i^2 + p_i^2 \le s^{\text {rated}}_i \},
\]
and the set $\mc U = \times_{i=1}^{18} \mc U_i$. Within normal operating ranges, there exists a function $y = \pi(u,w)$ that relates relates power injections (controllable and uncontrollable) to voltage magnitudes.\footnote{Treatment of robust stability of the OAG algorithm over a manifold $0 = \pi(u,w,y)$ is beyond the scope of this paper.} We formulate the optimization problem
\begin{align}\label{eq.opt.ps}
\begin{split}
\min_{u\in\mc U} ~& \|u-u_{\text{ref}}\|^2 + \sum_{i=1}^m \left\{ \max(0, \underline {y}-y_i, y_i-\overline {y})\right\}^2 \\
\mathrm{s.t.} ~& y = \pi(u,w).
\end{split}
\end{align}
where $u_{\text{ref}} \in \R^{36}$ is equal to the available PV power $p^{\max}_i$ for the element of $u$ corresponding to $p_i$ and $0$ for the element of $u$ corresponding to $q_i$, for all  $i=1,...,18$. 

\begin{figure}[htb]
\begin{center}
\input{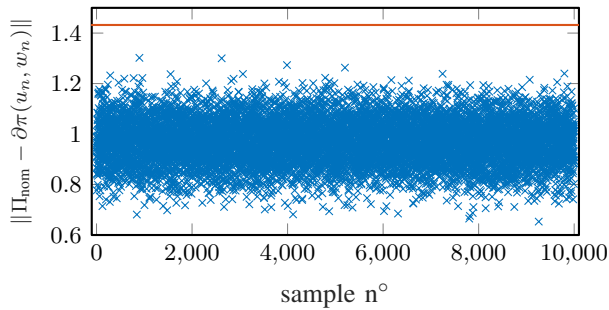}
\vspace{-0.3cm}
\caption{In order to estimate a norm bound $\gamma$ for the approximation error $\|\partial \pi(u,w) - \Pi_\text{nom}\|$ we sampled 10,000 operating points (by changing both $u$ and $w$). The red line in the plot show the chosen value of $\gamma=1.43$ used for the LMI test. Since the LMI test is feasible, we can guarantee robust stability for all $\pi$ such that $\|\partial \pi(u,w) - \Pi_\text{nom}\|\le \gamma$. }
\vspace{-0.3cm}
\label{fig.error}
\end{center}
\end{figure}

We set $\Pi_{\text{nom}}\in\R^{35\times 36}$ to be  Jacobian of the power flow equations of the feeder with the zero load profile and voltage magnitude of 1 p.u. at the PCC. In order to obtain a bound of the form~\eqref{eq.jacobian}, we sampled 10,000 operating points by randomly choosing both $u$ and $w$ from the power injections data and constructed the Jacobian from the controllable injections to the voltage magnitude using the method from~\cite{bernstein2018load}. The error $\|\partial \pi(u_n,w_n) - \Pi_\text{nom}\|$ of each sample is shown in Figure~\ref{fig.error}. As a safety factor, we multiplied the maximum empirical error observed by 1.1 to obtain $\gamma = 1.43$. 

In order to certify the stability of the OAG algorithm for all possible maps that satisfy~\eqref{eq.jacobian}, we follow the procedure outlined in Section~\ref{sec.application.optimization.lft}. In particular, we construct the matrices $(A,B,C,D)$ as in~\eqref{eq.abcd} and the cone $\mathbf \Theta$ as in~\eqref{eq.theta}. We use $\Pi = \Pi_{\text{nom}}$, and we can solve the single LMI~\eqref{eq.test} for $\Theta\in\mathbf \Theta$, with $P = I$ and $\rho = 0.45$ (0.59\,s using MOSEK on 2.5 GHz Intel Core i7 processor). 
%

\begin{figure}[htbp]
\begin{center}
\includegraphics[width=\columnwidth]{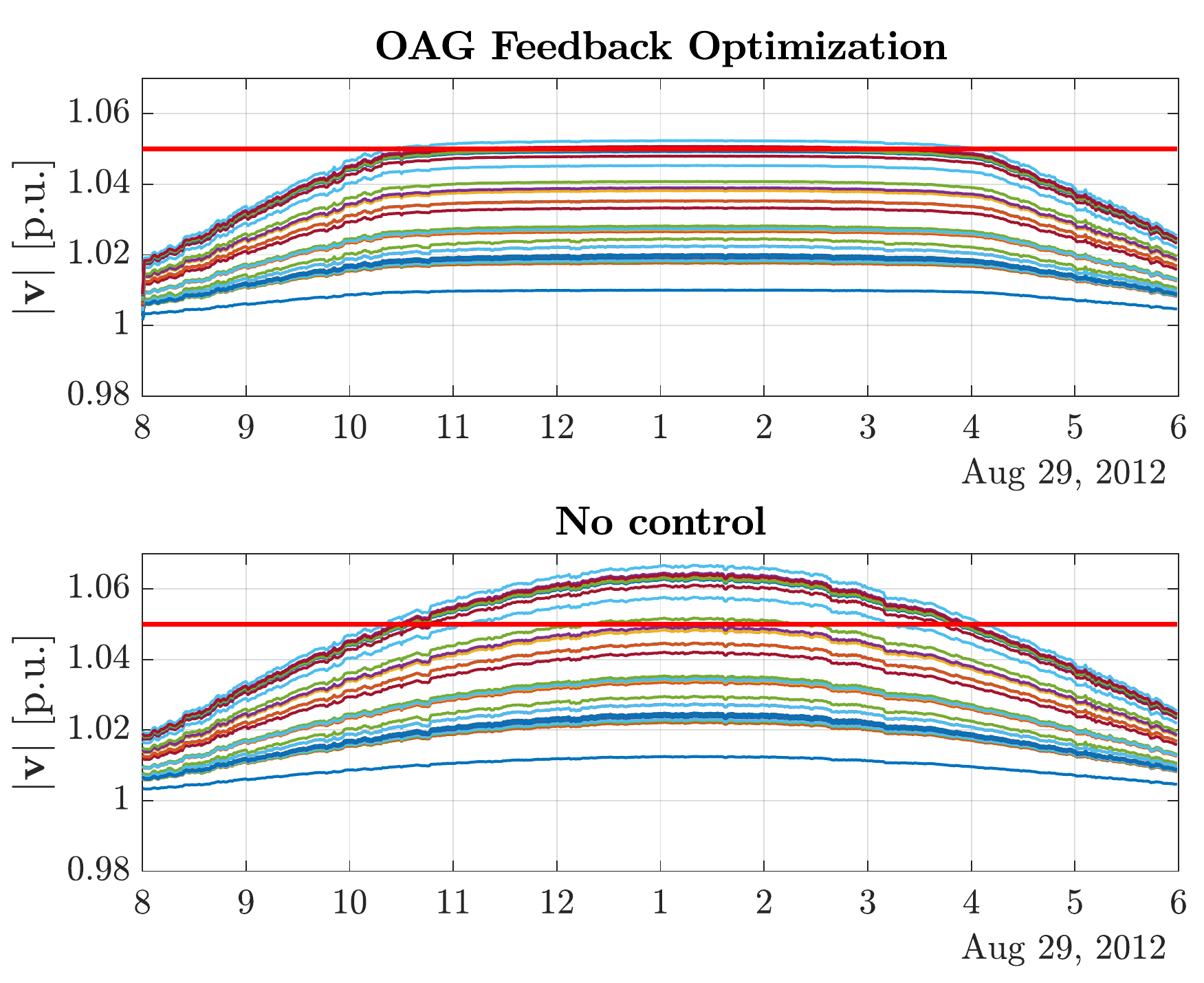}
\vspace{-0.6cm}
\caption{Comparison of  OAG vs no control for the IEEE37 test feeder. The disturbance $w$ (uncontrollable power injections and solar radiation) is taken from real data from Anatolia CA.}   
\vspace{-0.3cm}
\label{fig.comparison}
\end{center}
\end{figure}

By Proposition~\ref{propo.convergence}, the OAG algorithm is robustly stable with respect to the uncertainty and reaches the unique online approximate solution (which is, of course, different for every $w$). Figure~\ref{fig.comparison} shows a simulation of the OAG algorithm~for problem~\eqref{eq.opt.ps} applied to the IEEE37 bus system. The algorithm uses feedback on voltage measurements, which are computed by solving the AC power flow equations with MATPOWER at each time-step. The OAG is robustly stable (as predicted by the solvability of the LMI~\eqref{eq.test} and Proposition~\ref{propo.convergence}) and is able to significantly reduce over-voltage.

\section{Conclusions}
In this paper we studied a gradient-based optimization algorithm applied in feedback to a physical system. We characterized the equilibria of the feedback interconnection and we proposed a framework based on robust control theory to verify robust stability with respect to model mismatch and external disturbances. The results were illustrated on a realistic example from power systems. The first interesting extension is to use the LMI conditions from this paper to obtain linear approximations $\Pi$ which are robustly stable by design. Future work on will aslo focus on the robust stability analysis of more complex online algorithms that make better use of the available model information and real-time measurements. Further, for the specific case of power systems, we believe that combining tailored model uncertainty descriptions and pre-processing to identify (and remove) redundant constraints from the problem formulation (see e.g.~\cite{molzahn2019grid}) will certainly lead to less conservative robust stability guarantees.

\bibliographystyle{IEEEtran}
\bibliography{IEEEabrv,bib_file2}

\end{document}